\documentclass[12pt]{article}
\usepackage{amsmath,amssymb,amsthm,amsfonts}
\usepackage{color}
\usepackage{calligra}
\usepackage[T1]{fontenc}
\usepackage{graphicx}
 
\newcommand{\NN}{\mathbb N}
\newcommand{\RR}{{\mathbb R}}
\newcommand{\ZZ}{{\mathbb Z}}
\newcommand{\mA}{{\mathsf A}}
\newcommand{\mB}{{\mathsf B}}
\newcommand{\mC}{{\mathsf C}}
\newcommand{\mD}{{\mathsf D}}
\newcommand{\mF}{{\mathsf F}}
\newcommand{\mG}{{\mathsf G}}
\newcommand{\mR}{{\mathsf R}}
\newcommand{\mY}{{\mathsf Y}}
\newcommand{\mV}{{\mathsf V}}
\newcommand{\mZ}{{\mathsf Z}}

\newtheorem{theorem}{Theorem}[section]

\numberwithin{equation}{section}
\title{Quasi-Interpolant Operators and the Solution of Fractional Differential Problems}
\author{Enza Pellegrino\thanks{\noindent {\it Dip. DIIE, Universit\`{a} de L'Aquila }, Via G. Gronchi 18, 67100  L'Aquila, Italy. \rule{1cm}{0cm} \rule{0.5cm}{0cm} e-mail: \tt{enza.pellegrino@univaq.it}}, 
	Laura Pezza$^\dag$, Francesca Pitolli\thanks{\noindent {\it Dip. SBAI, Universit\`{a} di Roma ''La Sapienza''} Via A. Scarpa 16, 00161 Roma, Italy. 
		\rule{0.5cm}{0cm} e-mail: \tt{francesca.pitolli,laura.pezza@sbai.uniroma1.it} } }

\date{}

\begin{document}

\maketitle
	
\begin{abstract}
	Nowadays, fractional differential equations are a well established tool to model phenomena from the real world. Since the analytical solution is rarely available, there is a great effort in constructing efficient numerical methods for their solution. In this paper we are interested in solving boundary value problems having space derivative of fractional order. To this end, we present a collocation method in which the solution of the fractional problem is approximated by a spline quasi-interpolant operator. This allows us to construct the numerical solution in an easy way. We show through some numerical tests that the proposed method is efficient and accurate.
\end{abstract}

\section{Introduction}
In recent years fractional differential equations are becoming a powerful tool to describe real-world phenomena where nonlocality is a key ingredient. Starting from the fundamental paper  \cite{Ca67} by Caputo, where the fractional derivative was used for the first time to describe dissipation phenomena in Earth free modes, the literature on fractional models has exploded and now fractional differential equations are used in several fields, like continuum mechanics, signal processing, biophysics (see, \cite{KST06,Ma10,SKM93,Ta10} and references therein). The Caputo derivative is especially suitable to describe real phenomena since in many ways it behaves like the usual derivative of integer order. In particular, the Caputo derivative of constant functions is zero, which is not true for the Riemann-Liouville derivative \cite{Po99}. Moreover,  initial or boundary conditions can be easily applied \cite{Di10}. For details on fractional calculus see, for instance, \cite{Di10,OS74,Po99,SKM93}.
\\
Since the analytical solution of fractional differential equations can be rarely obtained explicitly, to solve these kinds of problems numerical methods are mandatory. There is a huge literature on numerical methods for fractional differential problems (see, \cite{Ba10,LC18,LZ15,Pi18a} and references therein). A crucial point to construct efficient methods is their ability to approximate the nonlocal behavior of the fractional derivative. In this respect, collocation methods that use information of the approximating function in all the discretization interval are received great attention in recent years \cite{KPT16,PT11,PP18.a,PP18.b,PP18c}.
\\
In this paper, we present a collocation method based on spline quasi-interpolant operators suitable to solve boundary value differential problems having fractional derivative in space. In particular, we are interested in solving linear boundary value problems of type
\begin{equation} \label{eq:fract_eq}
\left  \{ \begin{array}{ll}
D_x^{\gamma} y(x) + f(x)\,y(x) = g(x)\,, &0 < x < L\,,  \\ \\
\rho_{r0}\,y(0)+\rho_{r1}\,y'(0) + \zeta_{r0}\,y(L)+\zeta_{r1}\,y'(L) = c_r\,, & 1 \le r \le \lceil \gamma \rceil\,,
\end{array} \right.
\end{equation}
where $\gamma\in \bigl(\lfloor \gamma \rfloor,\lceil \gamma \rceil\bigr)$ is a given real number,
$f$ and $g$ are continuous given functions, and $\rho_{r0}$, $\rho_{r1}$, $\zeta_{r0}$, $\zeta_{r1}$, $c_r$ are given parameters. Here, we assume $L$ be a positive integer. Moreover, we assume the boundary conditions are linearly independent so that the differential problem has a unique solution \cite{Di10}.

The derivative appearing in the differential problem (\ref{eq:fract_eq}) should be intended in the Caputo sense. For a sufficiently smooth function, the Caputo fractional derivative is defined as 
\begin{equation} \label{eq:CDx}
D_x^\gamma \, y(x) := \frac1{\Gamma(\lceil \gamma \rceil-\gamma)} \, \int_0^x \frac{y^{(\lceil \gamma \rceil)}(\xi)}{(x-\xi)^{\gamma-\lceil \gamma \rceil+1}} \, d\xi\,, 
\end{equation} 
where $\Gamma$ is the Euler's gamma function
$$ 
\Gamma(\gamma) := \int_0^\infty \, \xi^{\gamma-1} \, {\rm e}^{-\xi} \, d\xi\,.
$$
In the method we present in this paper, we approximate the solution to the differential problem (\ref{eq:fract_eq}) by a spline quasi-interpolant. Polynomial spline quasi-interpolants are linear operators that are represented as a linear combination of spline basis functions whose coefficients are chosen in order to achieve some special properties, like shape preserving properties or good approximation order \cite{dBF73,LLM01,LS75,Sa05}. Thus, quasi-interpolants have a greater flexibility with respect to interpolation with the further advantage that they are easy to construct. 
We show that the proposed method is accurate and efficient since the fractional derivative of the approximating function can be evaluated explicitly. As a consequence, the nonlocal behavior of the fractional derivative can be easily taken into account.

The paper is organized as follows. The main properties of the B-spline basis we use to construct the approximating function are described in Section~\ref{sec:Bspline} while Section~\ref{sec:Bspline_der} is devoted to their fractional derivative. In Section~\ref{sec:QIop} we collect the main properties of the Schoenberg-Bernstein operator we use to approximate the solution of Equation~(\ref{eq:fract_eq}). The numerical method we propose is described in Section~\ref{sec:collmet} while some numerical results are shown in Section~\ref{sec:numtest}. Finally, some conclusions are drawn in Section~\ref{sec:concl}.

\section{The Cardinal B-splines}
\label{sec:Bspline}
	
The cardinal B-splines are piecewise polynomials with breakpoints at the integers \cite{dB78,Sc07}. They can be defined as
\begin{equation}  \label{B_spline}
B_n(x):=\frac{1}{n!}\,\Delta^{n+1} x_+^n, \qquad n \ge 0\,,
\end{equation}
where
\begin{equation} \label{eq:ffdo} 
\Delta^n \, f(x):= \sum_{\ell=0}^{n} \, (-1)^\ell \, {\dbinom{n}{\ell} \, f(x-\ell)}\,, \quad n \in \NN\,, 
\end{equation}
is the backward finite difference operator and $x_+^n :=  \bigl(\max(0,x) \bigr) ^n$ is the truncated power function. 
\\
The integer translates $\{B_n(x-\ell), \ell \in \ZZ\}$ form a basis for the spline space of degree $n$ on the whole line. A basis for the finite interval $[0,L]$, $L\ge n+1$, can be obtained by restriction, i.e.,
\begin{equation}
{\cal B}_n(x) = \{B_n(x-\ell), -n\le \ell \le L-1 \}\,, \qquad x \in [0,L]\,.
\end{equation}
We recall that the basis ${\cal B}_n(x)$ is totally positive, reproduces polynomial up to degree $n$ and is a partition of unity.

The B-spline bases can be generalized to any sequence of equidistant knots on the interval $[0,L]$ by mapping $x \to h^{-1}x$, where $h$ is the space step:
$$ 
{\cal B}_{h,n}(x) = \{B_{h,\ell,n}(x) = B_n(h^{-1}x -\ell), -n\le \ell \le h^{-1}L -1\}\,, \qquad x \in [0,L]\,.
$$
Thus, ${\cal B}_{h,n}$ is a basis for the spline space of degree $n$ having breakpoints at the knots $h\ell$, $0 \le \ell \le h^{-1}L$.	
We observe that for $-n\le \ell \le -1$ the functions $B_{h,\ell, n}(x)$ and $B_{h,h^{-1}L+\ell, n}(x)$ are left and right edge functions, respectively. Their support is $[0,h\ell]$ and $[L-h\ell,L]$, respectively. The functions $B_{h,\ell,n}(x)$ with $0\le \ell \le h^{-1}L-n-1$ are interior functions having support $[h\ell,h(\ell+n+1)]$.

\section{The Fractional Derivative of the Cardinal B-splines}
\label{sec:Bspline_der}
	
The fractional derivatives of the B-spline functions are fractional B-splines, i.e., piecewise polynomials of noninteger degree \cite{UB00}. 
Their explicit expression can be obtained by applying the Caputo differential operator (\ref{eq:CDx}) to the basis functions $B_{h,\ell,n}(x)$ (see \cite{PPP19,PP18.a} for details).

The Caputo derivative of the interior functions and of the right edge functions can be evaluated by the differentiation rule
\begin{equation} \label{eq:fract_der}
D^{\gamma}_x \, B_n (x)=  \frac{\Delta^{n+1} \, x_+^{n-\gamma}} {\Gamma(n+1-\gamma)}\,,  \qquad x\ge 0\,,  \qquad 0 < \gamma < n\,,
\end{equation}
where $x_+^\gamma = (\max(0,x) \bigr) ^\gamma$ is the fractional truncated power function (cf. \cite{PP18.a,UB00}). This formula generalizes to the noninteger case the well-known differentiation rule for the ordinary derivative of the B-spline
\begin{equation} \label{eq:int_der}
B_n^{(m)} (x)=  \frac{\Delta^{n+1} \, x_+^{n-m}} {(n-m)!}\,,  \qquad 0 \le m \le n-1\,.
\end{equation}
For  $-n \le \ell \le -1$, the Caputo derivative of the left edge functions is given by \cite{PPP19}
\begin{equation}\label{eq:edge_der}
\begin{array}{lcl}
\displaystyle D_{x}^{\gamma}B_{n,\ell}(x)&=&\displaystyle \frac{\Delta^{n+1}(x-\ell)_+^{n-\gamma}}{\Gamma(n+1-\gamma)} - \sum_{r=0}^{-\ell-1} (-1)^r \binom{n+1}{r}\bigg( \frac{(x-\ell-r)^{n-\gamma}}{\Gamma(n+1-\gamma)} - \\ \\
 \rule{1cm}{0cm}  &&\displaystyle \sum_{p=0}^{n-\lceil \gamma \rceil}  \frac{(-\ell-r)^{n-\lceil \gamma \rceil-p}\,x^{\lceil \gamma \rceil-\gamma+p}}{(n-\lceil \gamma \rceil-p)!\,\Gamma(\lceil \gamma \rceil-\gamma+p+1)} \bigg ).
\end{array}
\end{equation}
The fractional derivative of the refined basis functions $B_{h,\ell,n}$ can be evaluated recalling that, for any function $f$ sufficiently smooth, it holds
$$
D^{\gamma}_x f(h^{-1} x-\ell)=h^{-\gamma } D_{h^{-1} x}^{\gamma}f(h^{-1} x-\ell)
$$
(cf. \cite{PPP19}). 		

\section{Quasi-interpolant Operators}
\label{sec:QIop}

A quasi-interpolant operator is an approximation of a given function that reproduces polynomials up to a given degree.
In particular, a spline quasi-interpolant operator is a linear operator of type
\begin{equation}\label{eq:Q_n}
{\cal Q}_{n} \, y(x) = \sum_{\ell \in \ZZ} \, \mu_{\ell}(y) \, B_{n}(x-\ell)\,,
\end{equation}
where $\mu_{\ell}(y)$, $\ell \in \ZZ$, are continuous linear functionals that are determined by imposing that ${\cal Q}_{n} \, y$ is exact on polynomials up to degree $m \le n$. Usually, the functionals  $\mu_{\ell}(y)$ are assumed to be local, i.e.,
only values of $y$ in some neighborhood of $\sigma_{\ell,n}= supp \, B_n(x-\ell)$ are used to construct $\mu_{\ell}(y)$.
We notice that since $\mu_{\ell}(y)$ is local and $B_{n}(x-\ell)$ has compact support, for any $x\in \RR$ the sum in (\ref{eq:Q_n}) is actually a finite sum.

There are several kinds of quasi-interpolant spline operators (see, for instance, \cite{dBF73,GS88,LLM01,LS75,Sa05}). In this paper we consider Bernstein type operators \cite{Sa91} in which the functionals $\mu_{\ell}(y)$ are suitable values of $y$ evaluated on points belonging to $\sigma_{\ell,n}$. The simplest choice is
\begin{equation}
\mu_{\ell} (y) = y(\theta_{\ell})\,, 
\end{equation}
where
\begin{equation}
\theta_{\ell}=\ell + \frac {n+1}2\,, \qquad \ell \in \ZZ\,,
\end{equation}
are the Schoenberg nodes. This choice leads to the Schoenberg-Bernstein operator that reproduces linear functions and has approximation order 1 \cite{Sc67}. Even if the approximation order is poor, the Schoenberg-Bernstein operator has many properties useful in applications. In particular, it is a positive operator that has shape preserving properties. In fact, it enjoys the variation diminishing property, i.e., for any linear function $\Lambda$ and any function $y$ it holds
$$
S^-({\cal Q}_{n} \,( y-\Lambda)) \le S^-(y-\Lambda)\,,
$$
where $S^-(y)$ denotes the number of strict sign changes of the function $y$. This property reveals particular attractive in geometric modeling where the approximation of a given set of data is required to reproduce their shape \cite{Go96}.

The operator ${\cal Q}_{n} \, y$ is refinable, i.e., in the spline spaces generated by the B-spline basis ${\cal B}_{h,n}$, we can construct the refined operator
\begin{equation} \label{eq:Qhn}
{\cal Q}_{h,n} \, y(x) = \sum_{\ell \in \ZZ} \, \mu_{h,\ell}(y) \, B_{h,\ell,n}(x)\,,
\end{equation}
where $\mu_{h,\ell} \, (y)$ uses values of $y$ in $supp \, B_{h,\ell,n}$. The functionals $\mu_{h,\ell}(y)$ have expression
$$ 
\mu_{h,\ell} (y) = y(\theta_{h,\ell})\,,
$$
where $\theta_{h,\ell}= h\,\theta_\ell$ are the refined Schoenberg nodes.

\section{The Quasi-Interpolant Collocation Method}
\label{sec:collmet}

To solve the fractional differential problem (\ref{eq:fract_eq}) we approximate its solution by the refinable Schoenberg-Bernstein operator (\ref{eq:Qhn}) restricted to the interval $[0,L]$, i.e.,
\begin{equation}\label{eq:yapprox}
y(x)\approx y_{h,n}(x) = \sum_{\ell = -n}^{N_h} \, y_{h,n}(\tilde \theta_{h,\ell}) \, B_{h,\ell,n}(x)\,, \quad N_h = h^{-1}L-1\,, \qquad x \in [0,L]\,,
\end{equation}
where $\tilde \theta_{h,\ell}$ are the Schoenberg nodes for the interval $[0,L]$.
To determine the unknown coefficients $\{y_{h,n}(\tilde \theta_{h,\ell}),-n\le \ell \le N_h\}$ we solve the differential problem on a set of {\em  collocation points}. For the sake of simplicity, here we assume the collocation points are a set of equidistant nodes on the interval $[0,L]$ having distance $\delta=2^{-s}$,
\begin{equation}
X_\delta = \{x_{r}=\delta\,r,0\le r \le N_\delta\}\,, \qquad N_\delta = \delta^{-1}L\,.
\end{equation}
Thus, collocating Equation (\ref{eq:fract_eq}) on the nodes $X_\delta$ and using (\ref{eq:yapprox}) we get the linear system
\begin{equation} \label{eq:coll_eq}
\left  \{ \begin{array}{ll}
D_x^{\gamma} y_{h,n}(x_r) + f(x_r)\,y_{h,n}(x_r) = g(x_r)\,, \qquad 1 \le r \le N_\delta-1\,, \\ \\
\rho_{r0}\,y_{h,n}(x_0)+\rho_{r1}\,y'_{h,n}(x_0) + \zeta_{r0}\,y_{h,n}(x_{N_\delta})+\zeta_{r1}\,y'_{h,n}(x_{N_\delta}) = c_r\,, & 1 \le r \le \lceil \gamma \rceil\,.
\end{array} \right.
\end{equation}
Now, let 
\begin{equation*}
\mY_{h,\delta} = \bigl[y_{h,n}(\tilde \theta_{h,\ell}),-n\le \ell \le N_h \bigr]^T\,,
\end{equation*}
be  the unknown vector, 
\begin{equation*} \label{Ajs}
\mA_{h,\delta}= \bigl [B_{h,\ell,n}(x_r), 1 \le r\le N_\delta-1, -n \le \ell \le N_h \bigr]
\end{equation*}
and
\begin{equation*} \label{Djs}
\mD_{h,\delta}= \bigl [D^{\gamma}_x  B_{h,\ell,n}(x_r), 1 \le r\le N_\delta-1, -n \le \ell \le N_h \bigr]
\end{equation*}
be the collocation matrices of the refinable basis ${\cal B}_{h,n}$ and of its fractional derivative. Then, let 
\begin{equation*}
\mF_\delta = \bigl [f(x_r), 1 \le r\le N_\delta-1 \bigr]^T\,,
\qquad 
\mG_\delta = \bigl [g(x_r), 1 \le r\le N_\delta-1 \bigr]^T\,,
\end{equation*}
be the know terms. Finally, we define the vector parameters
\begin{equation*}
\begin{array}{ll}
\mR_k = \bigl [\rho_{rk},1\le r \le \lceil\gamma \rceil\bigr ]^T\,, & k=0,1\,, \\ \\
\mZ_k = \bigl [\zeta_{rk},1\le r \le \lceil\gamma \rceil\bigr ]^T\,, & k=0,1\,, \\ \\
\mC = \bigl [c_r,1\le r \le \lceil\gamma \rceil\bigr ]^T\,,
\end{array}
\end{equation*}
and the vectors containing the boundary values of the basis functions and of their first derivative
\begin{equation*} \label{Bjs}
\begin{array}{l}
\mB_{h,\delta}(x)=\bigl[B_{h,\ell,n}(x),-n \le \ell \le N_h \bigr]\,, \qquad x=0,L\,, \\ \\
\mB'_{h,\delta}(x)=\bigl[B'_{h,\ell,n}(x),-n \le \ell \le N_h \bigr]\,, \qquad x=0,L\,.
\end{array}
\end{equation*} 
Thus, Equation (\ref{eq:coll_eq}) can be written in matrix form as 
\begin{equation} \label{eq:matrix_coll}
\left \{
\begin{array}{ll}
(\mD_{h,\delta}+\mF_\delta\circ \mA_{h,\delta} )\, \mY_{h,\delta} = \mG_\delta\,,\\ \\
\bigl ( \mR_0\,\mB_{h,\delta}(0) + \mR_1\,\mB'_{h,\delta}(0) + \mZ_0\,\mB_{h,\delta}(L) + \mZ_1\,\mB'_{h,\delta}(L) \bigr )\, \mY_{h,\delta}=\mC\,,
\end{array} 
\right.   
\end{equation}
Here, $\mV\circ \mA$ denotes the entrywise product between a vector $\mV$ and a matrix $\mA$ meaning that $\mV$ has to be intended as a matrix having as many columns as $\mA$, each column being a replica of the vector $\mV$ itself.
The entries of the matrices $\mA_{h,\delta}$ and $\mD_{h,\delta}$ can be easily evaluated using formulas given in Sections~\ref{sec:Bspline}-\ref{sec:Bspline_der}.
\\
The linear system (\ref{eq:matrix_coll}) has $N_\delta-1+\lceil \gamma \rceil$ equations and $N_h+n+1$ unknowns. To guarantee the existence of a unique solution the refinement step $h$, the distance of the collocation points $\delta$ and the degree of the B-spline $n$ have to be chosen such that $N_\delta-1+\lceil \gamma \rceil \ge N_h+n+1$ \cite{PPP19}. We notice that the choice $N_\delta-1+\lceil \gamma \rceil > N_h+n+1$ is preferable since in this case there is a greater flexibility in the choice of the degree of the B-spline. In this case we get an overdetermined linear system that can be solved by the least squares method.

\noindent
Finally, following the same reasoning line as in \cite{PPP19} (cf. also \cite{As78}) it can be proved that the collocation method described above is convergent.

\begin{theorem}
	The collocation method is convergent, i.e. 
	$$
	\lim_{h \to 0}\| y(x) - y_{h,n}(x)\|_\infty = 0,
	$$
	where $\|y(x)\|_\infty = \max_{x \in [0,L]}|y(x)|$. 
\end{theorem}

\noindent
We notice that since the convergence order of spline collocation methods is related to the approximation properties of the spline spaces, we expect the infinity norm of the error to decrease at least as $h^\nu$, where $\nu$ is the smoothness of the known terms, providing that the approximating function and the differential operator are sufficiently smooth.

\section{Numerical Tests}
\label{sec:numtest}
		
\subsection{Example 1}

In the first test we solve the fractional differential problem
\begin{equation}
\left \{ \begin{array}{l}
D_x^\gamma y(x) + f(x) \, y(x) = g(x)\,, \qquad 0< x < 1 \,, \\ \\
y(0)+y(1) = 2\,,
\end{array} \right.
\end{equation}
where $\gamma \in (0,1)$ is a given real number and
$$
f(x)=x^{\frac 12}\,, \qquad g(x)=\frac{2}{\Gamma(2-\gamma)}\,x^{1-\gamma} + 2\,x^{\frac 32}\,.
$$
The exact solution is $y(x) = 2x$ so that the collocation method is exact for $n\ge 1$. To compare the exact and the numerical solutions, we evaluate the infinity norm of the error $e_{h,n}(x)=y(x)-y_{h,n}(x)$ as
$$
\|e_{h,n}\|_\infty = \max_{0\le r \le \eta N_\delta} |e_{h,n}(x_r)|\,,
$$
where $x_r = \delta r /\eta$, $0\le r \le \eta N_\delta$ with $\eta \in \NN^+$. In the tests we choose $\eta = 4$.
In the table below we list the infinity norm of the error we obtain using the Schoenberg-Bernstein operator with the B-splines of degree $n=3$ for $h=1/8$ and $\delta = 1/16$. To give an idea of the conditioning of the final linear system, the condition number $\kappa_{h,n}$ is also shown.
\\
\begin{center}
	\begin{tabular}{c|c|c|c}
		$\gamma$ & 0.25 & 0.5 & 0.75 \rule{0cm}{0.3cm}\\ 
		\hline  \rule{0cm}{0.5cm} 
		$\|e_{h,3}\|_\infty$ & 7.33e-15 & 1.09e-14 & 2.44e-15	\\
		$\kappa_{h,3}$       & 2.58e+01 & 1.67e+01 & 1.64e+01
	\end{tabular}
\end{center}		
As expected, the error is in the order of the machine precision.

\subsection{Example 2}

In the second test we solve the fractional differential problem
\begin{equation} \label{eq:fracdiffeq}
\left \{ \begin{array}{ll}
D_x^\gamma \, y(x) + y(x) = g(x) \,, & \qquad 0 < x < 1\,, \\ \\
y(0) = 0\,, \quad y(1) = 1\,, & 
\end{array} \right.
\end{equation}
where  $\gamma \in (1,2)$ is a given real number and
$$
g(x) = \frac{\Gamma(\nu+1)}{\Gamma(\nu+1-\gamma)}\,x^{\nu-\gamma}+x^\nu\,.
$$
The exact solution is $y(x) = x^\nu$. We approximate the solution by the Schoenberg-Bernstein operator with $n=4, 5, 6$ when $\nu = 2.5$ and $\gamma = 1.25, 1.5, 1.75$. The infinity norm of the error for different values of $h$ and $\delta = h/2$ is listed in the table below:
\begin{center}
	$\|e_{h,n}\|_\infty$ for $\gamma = 1.25$ \medskip\\ 
	\begin{tabular}{c|c|c|c|c}
		$h$ & $2^{-3}$ & $2^{-4}$ & $2^{-5}$ & $2^{-6}$ \rule{0cm}{0.3cm}\\ 
		\hline  \rule{0cm}{0.5cm} 
		$n=4$ & 5.91e-05 & 1.25e-05 & 2.67e-06 & 5.74e-07 \\
		\hline  \rule{0cm}{0.5cm} 
		$n=5$ & 3.84e-05 & 8.19e-06 & 1.73e-06 & 3.65e-07 \\	
		\hline  \rule{0cm}{0.5cm} 
		$n=6$ & 2.51e-05 & 5.35e-06 & 1.13e-06 & 2.39e-07	
	\end{tabular}
\end{center}
\begin{center}
$\|e_{h,n}\|_\infty$ for $\gamma = 1.5$ \medskip\\ 
	\begin{tabular}{c|c|c|c|c}
		$h$ & $2^{-3}$ & $2^{-4}$ & $2^{-5}$ & $2^{-6}$  \rule{0cm}{0.3cm}\\ 
		\hline  \rule{0cm}{0.5cm} 
		$n=4$ & 1.05e-04 & 2.68e-05 & 6.72e-06 & 1.67e-06 	\\
		\hline  \rule{0cm}{0.5cm} 
		$n=5$ & 6.68e-05 & 1.71e-05 & 4.33e-06 & 1.09e-06  \\	
		\hline  \rule{0cm}{0.5cm} 
		$n=6$ & 4.10e-05 & 1.05e-05 & 2.66e-06 & 6.68e-07 	
	\end{tabular}
\end{center}
\begin{center}
	$\|e_{h,n}\|_\infty$ for $\gamma = 1.75$ \medskip\\ 
	\begin{tabular}{c|c|c|c|c}
		$h$ & $2^{-3}$ & $2^{-4}$ & $2^{-5}$ & $2^{-6}$ \rule{0cm}{0.3cm}\\ 
		\hline  \rule{0cm}{0.5cm} 
		$n=4$ & 1.37e-04 & 4.20e-05 & 1.26e-05 & 3.72e-06 	\\
		\hline  \rule{0cm}{0.5cm} 
		$n=5$ & 8.02e-05 & 2.50e-05 & 7.62e-06 & 2.29e-06  \\	
		\hline  \rule{0cm}{0.5cm} 
		$n=6$ & 4.34e-05 & 1.34e-05 & 4.06e-06 & 1.22e-06 	
	\end{tabular}
\end{center}
As expected, the norm of the error decreases when $h$ decreases. We notice that the error decreases also when $n$ increases.

\subsection{Example 3} 

In the last test we solve the fractional differential problem
\begin{equation} \label{eq:fracdiffeq}
\left \{ \begin{array}{ll}
D_x^\gamma \, y(x) + y(x) = 0 \,, & \qquad 0 < x < 1\,, \\ \\
y(0) = 1\,, \quad y(1) = E_\gamma(-1^\gamma)\,, & 
\end{array} \right.
\end{equation}
where $\gamma \in (1,2)$ is a given real number and
$$
E_\gamma(x) = \sum_{\ell\ge 0} \frac{x^\ell}{\Gamma(\gamma \ell +1)}\,, 
$$
is the one-parameter Mittag-Leffler function \cite{GKMR14}.
The exact solution is $y(x) =  E_\gamma(-x^\gamma)$. We approximate the solution by the Schoenberg-Bernstein operator with $n=4, 5, 6$. The infinity norm of the error for different values of $h$ and $\delta = h/2$  when $\gamma = 1.25, 1.5, 1.75$ is listed in the tables below:
\begin{center}
$\|e_{h,n}\|_\infty$ for $\gamma = 1.25$ \medskip\\ 
	\begin{tabular}{c|c|c|c|c}
		$h$ & $2^{-3}$ & $2^{-4}$ & $2^{-5}$ & $2^{-6}$  \\ 
		\hline  \rule{0cm}{0.5cm} 
		$n=4$ & 9.01e-03 & 4.55e-03 & 2.27e-03 & 1.13e-03 	\\
		\hline  \rule{0cm}{0.5cm} 
		$n=5$ & 8.01e-03 & 4.05e-03 & 2.03e-03 & 1.02e-03 \\	
		\hline  \rule{0cm}{0.5cm} 
		$n=6$ & 7.06e-03 & 3.57e-03 & 1.79e-03 & 8.98e-04	
	\end{tabular}
\end{center}
\begin{center}
	$\|e_{h,n}\|_\infty$ for $\gamma = 1.5$ \medskip\\
	\begin{tabular}{c|c|c|c|c}
		$h$ & $2^{-3}$ & $2^{-4}$ & $2^{-5}$ & $2^{-6}$  \rule{0cm}{0.3cm}\\ 
		\hline  \rule{0cm}{0.5cm} 
		$n=4$ & 3.66e-03 & 1.87e-03 & 9.41e-04 & 4.72e-04 	\\
		\hline  \rule{0cm}{0.5cm} 
		$n=5$ & 3.11e-03 & 1.59e-03 & 8.02e-04 & 4.03e-04	\\
		\hline  \rule{0cm}{0.5cm} 
		$n=6$ & 2.63e-03 & 1.34e-03 & 6.76e-04 & 3.40e-04	
	\end{tabular}
\end{center}
\begin{center}
	$\|e_{h,n}\|_\infty$ for $\gamma = 1.75$ \medskip\\
	\begin{tabular}{c|c|c|c|c}
		$h$ & $2^{-3}$ & $2^{-4}$ & $2^{-5}$ & $2^{-6}$  \rule{0cm}{0.3cm}\\ 
		\hline  \rule{0cm}{0.5cm} 
		$n=4$ & 7.91e-04 & 4.10e-04 & 2.08e-04 & 1.05e-04 	\\
		\hline  \rule{0cm}{0.5cm} 
		$n=5$ & 6.32e-04 & 3.28e-04 & 1.67e-04 & 8.44e-05   \\	
		\hline  \rule{0cm}{0.5cm} 
		$n=6$ & 4.98e-04 & 2.57e-04 & 1.30e-04 & 6.57e-05	
	\end{tabular}
\end{center}
Also in this case the norm of the error decreases when $h$ decreases and $n$ increases. 

\section{Conclusion}
\label{sec:concl}
We have presented a collocation method based on spline quasi-interpolant operators. The method is easy to implement and has proved to be convergent. The numerical tests show that it is efficient and accurate. The method can be improved in several ways. First of all, we used the truncated B-spline bases to construct the Schoenberg-Bernstein operator. It is well known that truncated bases have low accuracy in approximating the boundary conditions which can result in a poor approximation of the solution. This problem can be overcome by using B-spline bases with multiple nodes at the endpoints of the interval. The use of this kind of B-splines for the solution of fractional problems has already considered in \cite{Pi18b} where the analytical expression of their fractional derivative is also given. As for the quasi-interpolants, even if the Schoenberg-Bernstein operator produces good results, it is just second order accurate. To increase the approximation order, different quasi-interpolant operators can be used, like the projector quasi-interpolants introduced in \cite{LLM01} or integral or discrete operators \cite{LS75,Sa05}. 
Finally, we notice that the accuracy of the method could also be improved by using Gaussian points instead of equidistant points (cf. \cite{As78,dBS73}).
These issues are at present under study.

%
%

\end{document}